# Distributive properties of the rationals

*Konstantine Zelator*
*Department of Mathematics and Computer Science*
*Rhode Island College*
*600 Mount Pleasant Avenue*
*Providence, R.I. 02908-1991, U.S.A.*
*e-mail address: 1) Kzelator@ric.edu*
*2) Konstantine_zelator@yahoo.com*

*Effective August 1, 2009, and for the academic year 2009-10;*
*Konstantine Zelator*
*Department of Mathematics*
*301 Thackeray Hall*
*139 University Place*
*University of Pittsburgh*
*Pittsburgh, PA 15260*
*U.S.A*
*e-mail address : kzet159@pitt.edu*

## 1. Introduction

Well before the intricacies of irrational numbers come into play in high school math education, a student learns about the fundamental properties of the four basic operations in the number system of the rational numbers. These basic binary operations are addition +, multiplication $\cdot$ , subtraction − , and division ÷. Two of the fundamental properties are the distributive properties of multiplication over addition and subtraction:

For any rational numbers $r_1, r_2, r_3$;

$$r_1 \cdot (r_2 + r_3) = r_1 \cdot r_2 + r_1 \cdot r_3 \quad (1)$$

$$\text{and} \quad r_1 \cdot (r_2 - r_3) = r_1 \cdot r_2 - r_1 \cdot r_3 \quad (2)$$

In a more general setting, if $*$ and $\circ$ are two binary operations in a number system; then the number system is said to have the distributive property of the operation $\circ$ over the operation $*$ if for any three elements $a, b, c$ in the number system,

$$a \circ (b * c) = (a \circ b) * (a \circ c). \quad (3)$$

In the case of the rationals, addition, multiplication and subtraction are defined for any ordered pair of rational numbers; while division will produce a rational number outcome, only when the second number ( in the ordered pair) is nonzero.

Now, consider the following question. In the rational number system, if $*$ is one of the four operations $+, \cdot, -, \div;$ and likewise, if $\circ$ is also one of these four operations. When is (3) satisfied, that is, for which rational triples is (3) satisfied? If we allow $\circ$ and $*$ to be the same operation, then there are preuisely 16 combinations with $*$ and $\circ$ being among the four basic operations. The two cases in which $\circ$ is the multiplication operation and $*$ is either addition or subtraction , are already answered by (1) and (2). Of the remaining 14 cases; 11 are easy/straight forward in that, almost always, one of the three rationals involved, must equal zero or 1. The last three cases are really interesting, and number theory is involved.

It appears that the determination of all such rational triples in each of these three cases is a very complicated, involved process. Instead, what we do in these three cases, is finding infinite families of rational triples with the said distributive property.

This is then the aim of this work.

To determine, fully or partly, the triples of rational numbers which satisfy the distributive property, for each of the fourteen combinations or cases.

## 2. Eleven easy/straight forward cases

1. *Distributive property of addition over itself*

    In this case, $r_1 + (r_2 + r_3) = (r_1 + r_2) + (r_1 + r_3) \Leftrightarrow r_1 = 0$

    *Rational triples:* $(0, r_2, r_3)$

2. *Distributive property of addition over subtraction*

    We have, $r_1 + (r_2 - r_3) = (r_1 + r_2) - (r_1 + r_3) \Leftrightarrow r_1 = 0$

    *Rational triples :* $(0, r_2, r_3)$

3. *Distributive property of multiplication over itself*

    This is, $r_1 \cdot (r_2 \cdot r_3) = (r_1 \cdot r_2) \cdot (r_1 \cdot r_3) \Leftrightarrow r_1 r_2 r_3 (r_1 - 1) = 0 \Leftrightarrow$
    $\Leftrightarrow (r_1 r_2 r_3 = 0 \text{ or } r_1 = 1)$. We have two nondisjoint families:

    *Family 1: rational triples* $(r_1, r_2, r_3)$ *with* $r_1 r_2 r_3 = 0$

    *Family 2: rational triples* $(1, r_2, r_3)$

4. *Distributive property of multiplication over division*

We have , $r_1 \cdot (r_2 \div r_3) = (r_1 \cdot r_2) \div (r_1 \cdot r_3)$,
with $r_1 r_3 \neq 0$.

Equivalently , $\left( \dfrac{r_1 r_2}{r_3} = \dfrac{r_1 r_2}{r_1 r_3}, r_1 r_3 \neq 0 \right) \Leftrightarrow \left\{ \begin{array}{l} r_2 (r_1 - 1) = 0 \\ r_1 r_3 \neq 0 \end{array} \right\} \Leftrightarrow$
$\Leftrightarrow ((r_2 = 0 \text{ or } r_1 = 1) \text{ and } r_1 r_3 \neq 0)$.

Rational triples : Two nondisjoint families:

*Family 1: Rational triples* $(r_1, 0, r_3)$ *with* $r_1 r_3 \neq 0$
*Family 2: Rational triples* $(1, r_2, r_3)$ *with* $r_3 \neq 0$

5. *Distributive property of subtraction over itself*

This is the case of $r_1 - (r_2 - r_3) = (r_1 - r_2) - (r_1 - r_3) \Leftrightarrow r_1 = 0$

*Rational triples :* $(0, r_2, r_3)$

6. *Distributive property of subtraction over addition*

In this case, $r_1 - (r_2 + r_3) = (r_1 - r_2) + (r_1 - r_3) \Leftrightarrow r_1 = 0$

*Rational triples :* $(0, r_2, r_3)$

7. Distributive property of division over itself

This property says that $r_1 \div (r_2 \div r_3) = (r_1 \div r_2) \div (r_1 \div r_3)$,

with $r_2 r_3 \neq 0$. This statement is equivalent to,

$\left( \dfrac{r_1 r_3}{r_2} = \dfrac{r_3}{r_2} \text{ and } r_1 r_2 r_3 \neq 0 \right) \Leftrightarrow (r_1 = 1 \text{ and } r_2 r_3 \neq 0)$.

*Rational triples:* $(1, r_2, r_3)$ with $r_2 r_3 \neq 0$.

8. Distributive property of division over multiplication.

This property says that $r_1 \div (r_2 \cdot r_3) = (r_1 \div r_2) \cdot (r_1 \div r_3)$, with $r_2 r_3 \neq 0$

Equivalently, $\left( \dfrac{r_1}{r_2 r_3} = \dfrac{r_1^2}{r_2 r_3} \text{ and } r_2 r_3 \neq 0 \right) \Leftrightarrow (r_1(r_1 - 1) = 0 \text{ and } r_2 r_3 \neq 0)$

$\Leftrightarrow ((r_1 = 0 \text{ or } r_1 = 1) \text{ and } r_2 r_3 \neq 0)$, We obtain two disjoint families:

*Family 1 : Rational triples* $(0, r_2, r_3)$ with $r_2 r_3 \neq 0$

*Family 2 : Rational triples* $(1, r_2, r_3)$ with $r_2 r_3 \neq 0$

9. Distributive property of division over addition

We have, $r_1 \div (r_2 + r_3) = (r_1 \div r_2) + (r_1 \div r_3)$, with $r_2 r_3 \neq 0$ and $r_2 \neq -r_3$.

Equivalently, $\left( \dfrac{r_1}{r_2 + r_3} = \dfrac{r_1}{r_2} + \dfrac{r_1}{r_3}, \ r_2 r_3 \neq 0, r_2 \neq -r_3 \right)$.

One possibility is $r_1 = 0$. Now, if $r_1 \neq 0$, the last equation is equivalent to

$r_2 r_3 = r_2 + r_3 \Leftrightarrow (divide \ with \ r_2 r_3) \ \dfrac{r_2}{r_3} + \dfrac{r_3}{r_2} - 1 = 0$, which shows that the rational number $\dfrac{r_2}{r_3}$ is a root of the quadratic equation $x^2 - x + 1 = 0$. But this quadratic equation has in fact, no real roots.

Thus, $r_1 = 0$.

*Rational triples :* $(0, r_2, r_3)$ with $r_2 r_3 \neq 0$

### 10. Distributive property of division over subtraction

In this situation, $r_1 \div (r_2 - r_3) = (r_1 \div r_2) - (r_1 \div r_3)$, with $r_2 r_3 \neq 0$ and $r_2 \neq r_3$. We have,

$$\frac{r_1}{r_2 - r_3} = \frac{r_1}{r_2} - \frac{r_1}{r_3}$$

One possibility is $r_1 = 0$. If $r_1 \neq 0$, then the last equation is equivalent to
$(r_2 r_3 = (r_3 - r_2)(r_2 - r_3); \text{ and with } r_2 r_3 \neq 0, r_2 \neq r_3)$

Equivalently, $r_2^2 - r_2 r_3 + r_3^2 = 0$; $\left(\frac{r_2}{r_3}\right)^2 - \left(\frac{r_2}{r_3}\right) + 1 = 0$, an impossibility since the quadratic equation $x^2 - x + 1 = 0$ has no real roots.

*Rational triples:* $(0, r_2, r_3)$ with $r_2 r_3 \neq 0$

### 11. Distributive property of addition over multiplication

This property says that
$$r_1 + (r_2 \cdot r_3) = (r_1 + r_2) \cdot (r_1 + r_3) \Leftrightarrow$$
$$\Leftrightarrow r_1 + r_2 r_3 = r_1^2 + r_1 r_3 + r_2 r_1 + r_2 r_3 \Leftrightarrow$$
$$\Leftrightarrow r_1 (r_1 + r_2 + r_3 - 1) = 0 \Leftrightarrow (r_1 = 0 \text{ or } r_1 + r_2 + r_3 = 1)$$

We obtain two nondisjoint families.

*Family 1: Rational triples* $(0, r_2, r_3)$

*Family 2: Rational triples* $(1 - (r_2 + r_3), r_2, r_3)$

## 3. Three interesting, complicated combinations

### 12. Distributive property of subtraction over multiplication

$$r_1 - (r_2 \cdot r_3) = (r_1 - r_2) \cdot (r_1 - r_3).$$

Equivalently, after some algebra,

$$r_1^2 - r_1 r_3 - r_1 r_2 + 2 r_2 r_3 - r_1 = 0 \tag{4}$$

First, we determine those triples $(r_1, r_2, r_3)$ which satisfy $r_1 r_2 r_3 = 0$ and (4):

If $r_1 = 0$, then also $r_2 r_3 = 0$. If $r_1 \neq 0$, and $r_2 = 0$ then (4) is equivalent to

$r_1(r_1 - r_3 - 1) = 0 \Leftrightarrow (\text{since } r_1 \neq 0) \; r_1 = r_3 + 1; \; \text{with } r_3 \neq -1 \; (\text{since } r_1 \neq 0)$

If $r_1 \neq 0$ and $r_3 = 0$, (4) similarly yields $r_1 = r_2 + 1.$, with $r_2 \neq -1$.

*Family 1: Rational triples* $(0, r_2, 0)$ *and* $(0, 0, r_3)$

*Family 2: Rational triples* $(r_3 + 1, 0, r_3)$ *with* $r_3 \neq -1$

*Family 3: Rational triples* $(r_2 + 1, r_2, 0)$, *with* $r_2 \neq -1$.

Next, assume $r_1 r_2 r_3 \neq 0$; each of $r_1, r_2, r_3$ is nonzero. Every nonzero rational can be uniquely written in the form $\dfrac{n}{d}$, where $d$ is a positive integer, $n$ is an integer, and $(n, d) = 1$ (i.e. $n$ and $d$ are relatively prime; their greatest common divisor being 1). Accordingly we set,

$$\left\{ \begin{array}{l} r_1 = \dfrac{n_1}{d_1}, r_2 = \dfrac{n_2}{d_2}, r_3 = \dfrac{n_3}{d_3}, \text{where } n_1, n_2, n_3 \text{ are nonzero integers,} \\ d_1, d_2, d_3 \text{ are positive integers and } (n_1, d_1) = (n_2, d_2) = (n_3, d_3) = 1 \end{array} \right\} \quad (5)$$

By (4) and (5) we obtain,

$$n_1^2 d_2 d_3 - n_1 n_3 d_1 d_2 - n_1 n_2 d_1 d_3 + 2 n_2 n_3 d_1^2 - n_1 d_1 d_2 d_3 = 0 \quad (6)$$

To determine all the positive integer solutions in $n_1, n_2, n_3, d_1, d_2, d_3$ of equation (6) and which satisfy conditions (5); is well beyond the scope of this work. However, we can find some families of solutions with $d_1 = d_2 = d_3 = 1$ (so all three rationals are nonzero integers).

Accordingly, with $d_1 = d_2 = d_3 = 1$, (6) takes the form

$$n_1^2 - n_1 n_3 - n_1 n_2 + 2 n_2 n_3 - n_1 = 0 \quad (7)$$

Let $\delta$ be the greatest common divisor of $n_1$ and $n_2$: $\delta = (n_1, n_2)$.

We have $n_1 = \delta \cdot N_1$, $n_2 = \delta \cdot N_2$; where $N_1$ and $N_2$ are relatively prime integers; $(N_1, N_2) = 1$.
Altogether,

$$n_1 = \delta \cdot N_1, n_2 = \delta \cdot N_2, (N_1, N_2) = 1 \quad (8)$$

From (7) and (8) we obtain,

$$\delta \cdot N_1^2 - N_1 \cdot n_3 - \delta \cdot N_1 N_2 + 2 N_2 \cdot n_3 - N_1 = 0;$$
$$\text{or equivalently} \quad N_1 \cdot (-\delta N_1 + n_3 + \delta N_2 + 1) = 2 N_2 n_3 \quad (9)$$

Equation (9) clearly shows that $N_1$ is a divisor of the product $2N_2 n_3$. Thus if we assume $N_1$ to be odd; then in view of $(N_1, N_2) = 1$, it follows that $(N_1, 2N_2) = 1$. Since $N_1$ is relatively prime to $2N_2$ and it divides the product $(2N_2) \cdot n_3$; it follows by Euclid Lemma, that $N_1$ must be a divisor of $n_3$. Euclid Lemma can be found in number theory texts, for example in **[1].**

Put $n_3 = N_1 \cdot N_3$, for some integer $N_3$ (10)

By (9) and (10) we get

$$-\delta N_1 + N_1 N_3 + \delta N_2 + 1 = 2N_2 N_3; \text{ or equivalently,}$$

$$\delta(N_1 - N_2) + N_3(2N_2 - N_1) = 1 \quad (11)$$

Note that since $(N_1, N_2) = 1$, we must have $(N_1 - N_2, 2N_2 - N_1) = 1$. It is now clear that for given relatively prime non zero integers $N_1$ and $N_2$, and with $N_1$ odd; equation (11) can be viewed as a linear diophantine equation in the variables $\delta$ and $N_3$. Such equations can be solved with a standard well known procedure or method, whose description can be found in number theory books ( see **[1]**). In the case of (11), we also have the restriction $\delta \geq 1$ (since $\delta$ is a greatest common divisor).

If we take $N_1 = 3$ and $N_2 = 2$, then from (11) we obtain

$$\delta + N_3 = 1; \quad N_3 = 1 - \delta$$

With $N_1 = 3, N_2 = 2, \delta \geq 1, N_3 = 1 - \delta$; and thus from (10) and (8),

$$n_3 = 3(1-\delta), \; n_1 = 3\delta, \; n_2 = 2\delta.$$

And form $d_1 = d_2 = d_3 = 1$ and (5), we obtain the following.

*Family 4: Rational triples*: $(r_1, r_2, r_3) = (3\delta, 2\delta, 3(1-\delta))$; *with $\delta$ being a positive integer greater that 1*

$$(\text{since } r_3 \neq 0)$$

*13. Distributive property of addition over division*

This is property states that $r_1 + (r_2 \div r_3) = (r_1 + r_2) \div (r_1 + r_3)$, with $r_3 \neq 0$, and $r_1 \neq -r_3$.

We have $r_1 + \dfrac{r_2}{r_3} = \dfrac{r_1 + r_2}{r_1 + r_3}$. First we treat the cases $r_1 r_2 r_3 = 0$, which means, since $r_3 \neq 0$, that $r_1 r_2 = 0$. If $r_1 = 0$, we obtain $\dfrac{r_2}{r_3} = \dfrac{r_2}{r_3}$, which is true. If $r_1 \neq 0$, then $r_2 = 0$, and so

$r_1 = \dfrac{r_1}{r_1 + r_3}$, which implies $(\text{since } r_1 \neq 0)$, $r_1 + r_3 = 1$; $r_1 = 1 - r_3$, with $r_3 \neq 1$.

We obtain the following families:

*Family 1:* Rational triples $(r_1, r_2, r_3) = (0, r_2, r_3)$ with $r_3 \neq 0$.

*Family 2:* Rational triples $(r_1, r_2, r_3) = (1 - r_3, 0, r_3)$ with $r_3 \neq 0, 1$.

Next, assume $r_1 r_2 r_3 \neq 0$, and set $r_1 = \dfrac{a}{b}, r_2 = \dfrac{c}{d}, r_3 = \dfrac{e}{f}$, where the denominators $b, d, f$ are positive integers and the numerators $a, c, e$ are nonzero integers relatively prime to $b, d$ and $f$ respectively. We go back to,

$$r_1 + \frac{r_2}{r_3} = \frac{r_1 + r_2}{r_1 + r_3} \tag{12}$$

Equivalently, $r_3 (r_1 + r_3) \left( r_1 + \dfrac{r_2}{r_3} \right) = r_3 (r_1 + r_3) \left( \dfrac{r_1 + r_2}{r_1 + r_3} \right)$; with $r_1 r_2 r_3 \neq 0$ and $r_1 \neq -r_3$.

$$r_3 r_1 (r_1 + r_3) + r_2 (r_1 + r_3) = r_3 (r_1 + r_2);$$
$$r_3 r_1^2 + r_1 r_3^2 + r_2 r_1 + r_2 r_3 = r_3 r_1 + r_3 r_2;$$
$$r_3 r_1^2 + r_1 r_3^2 + r_2 r_1 - r_3 r_1 = 0;$$
$$r_1 (r_1 r_3 + r_3^2 + r_2 - r_3) = 0; \text{ and since } r_1 \neq 0, \text{ we obtain}$$

$$r_1 r_3 + r_3^2 + r_2 - r_3 = 0 \tag{13}$$

Substituting for $r_1 = \dfrac{a}{b}, r_2 = \dfrac{c}{d},$ and $r_3 = \dfrac{e}{f}$ yields,

$$\begin{cases} adef + bde^2 + bcf^2 - bdef = 0 \\ \text{with } a, c, e \in \mathbb{Z}; b, d, f \in \mathbb{Z}^+, \\ \text{and } (a, b) = 1, (c, d) = 1, (e, f) = 1; \\ \text{and also } ace \neq 0 \end{cases} \tag{14}$$

If we take $b = e = f = 1$ in (14). Then,
$$ad + d + c - d = 0; \quad c = -ad \quad \text{we obtain the following family:}$$

*Family 3:* Rational triples $(r_1, r_2, r_3) = (a, -a, 1)$, with $a$ being a nonzero integer.

Going back to (14) and taking $e = f = 1$; we get
$$ad + bd + bc - bd = 0 \Leftrightarrow ad = -bc; \text{ and } abcd \neq 0;$$

we obtain $\dfrac{a}{b} = -\dfrac{c}{d}$. We have,

*Family 4: Rational triples* $(r_1, r_2, r_3) = \left(\dfrac{c}{d}, -\dfrac{c}{d}, 1\right)$, *with d a positive integer and c a nonzero integer*

Of course, Family 3 is a subset of Family 4.

Next, we go back to (14) and we take $d = b = 1$. We obtain quadratic equation in *e* with leading coefficient 1:

$$e^2 + f(a-1)e + cf^2 = 0 \tag{15}$$

The discriminant must be a perfect square:

$$[f(a-1)]^2 - 4cf^2 = K^2 \tag{16}$$

for some nonnegative integer $K$.

Accordingly (quadratic formula), $e = \dfrac{-f(a-1) \pm K}{2}$ \hfill (17)

From (16), $c = \dfrac{[f(a-1)]^2 - K^2}{4f^2}$

The constraints, according to (14), are $af \neq 0$ and, $-f(a-1) \pm K \neq 0$.

*Family 5: Rational triples* $(r_1, r_2, r_3) = \left(a, \dfrac{[f(a-1)]^2 - K^2}{4f^2}, \dfrac{-f(a-1) \pm K}{2f}\right)$

### 14. Distributive of subtraction over division

The stated property is $r_1 - (r_2 \div r_3) = (r_1 - r_2) \div (r_1 - r_3)$, with $r_3 \neq 0$ and $r_1 \neq r_3$.

If $r_1 r_2 r_3 = 0$, then $r_1 = 0$ or $r_2 = 0$. If $r_1 = 0$, we obtain

$-\dfrac{r_2}{r_3} = -\dfrac{r_2}{r_3}$, which is true.

If $r_1 \neq 0$, then $r_2 = 0$, and so $r_1 = \dfrac{r_1}{r_1 - r_3}$; and since $r_1 \neq 0$, we get $r_1 - r_3 = 1$, $r_1 = r_3 + 1$. We obtain the families,

*Family 1: Rational triples* $(r_1, r_2, r_3) = (0, r_2, r_3)$, *with* $r_3 \neq 0$.
*Family 2: Rational triples* $(r_1, r_2, r_3) = (r_3 + 1, 0, r_3)$, *with* $r_3 \neq 0, -1$.

Next, assume that $r_1 r_2 r_3 \neq 0$ and $r_1 \neq r_3$. We write the above property in the form,

$$r_1 - \frac{r_2}{r_3} = \frac{r_1 - r_2}{r_1 - r_3}.$$

We multiply the last equation by $r_3(r_1 - r_3)$ :

$$r_1 r_3 (r_1 - r_3) - r_2 (r_1 - r_3) = r_3 (r_1 - r_2);$$
$$r_1^2 r_3 - r_1 r_3^2 - r_2 r_1 + r_2 r_3 = r_3 r_1 - r_3 r_2;$$
$$r_1 r_3^2 + r_1 r_2 - 2 r_2 r_3 + r_1 r_3 - r_3 r_1^2 = 0$$

We put $r_1 = \frac{A}{B}$, $r_2 = \frac{C}{D}$, $r_3 = \frac{E}{F}$, with these fractions being in lowest terms and the denominators being positive. We have,

$$\frac{A}{B} \cdot \frac{E^2}{F^2} + \frac{A}{B} \cdot \frac{C}{D} - 2\frac{C}{D} \cdot \frac{E}{F} + \frac{A}{B} \cdot \frac{E}{F} - \frac{E}{F} \cdot \frac{A^2}{B^2} = 0$$

Multiply the last equation by $DB^2 F^2$ : we obtain,

$$\begin{cases} ABDE^2 + ABCF^2 - 2B^2 CEF + ABDEF - A^2 DEF = 0 \\ \text{with } B, D, F \text{ positive intergers}; A, C, E \text{ nonzero integers} \\ \text{and } (A, B) = (C, D) = (E, F) = 1 \end{cases} \quad (18)$$

If we take $A = B = C = 1$ in (18), then

$$DE^2 + F^2 - 2EF + DEF - DEF = 0;$$
$$F^2 - 2EF + DE^2 = 0; \text{ or equivalently, solving for } D \text{ produces}$$
$$D = \frac{F^2 - 2EF}{E^2} = \frac{F(F - 2E)}{E^2}.$$

Therefore we have the following family of rational triples.

*Family 3: Rational triples* $(r_1, r_2, r_3) = \left(1, \frac{E^2}{F(F - 2E)}, \frac{E}{F}\right)$, with $F$ being a positive integer, $E$ a nonzero integer, and *such that* $F > 2E$, *and* $(E, F) = 1$

# References

**[1]** Kenneth H.Rosen, *Elementary Number Theory and its Applications,* Fifth Edition, 721 p.p., Pearson , Addison Wesley 2005 ISBN 0-321-23707-2

For Euclid's Lemma, see Lemma 3.4 on page 109

For two-variable linear Diophantine equations, see Theorem 3.23 on page 134